\documentstyle[amssymb,amsfonts]{amsart}

\renewcommand{\a}{\alpha}
\renewcommand{\b}{\beta}
\newcommand{\ga}{\gamma}
\newcommand{\de}{\delta}

\newcommand{\varep}{\varepsilon}

\newcommand{\si}{\sigma}
\newcommand{\om}{\omega}

\newcommand{\Ga}{\Gamma}
\newcommand{\De}{\Delta}

\newcommand{\pr}{\partial}

\newcommand{\rr}{\mbox{{\bf R}}}  
\newcommand{\Z}{\mbox{{\bf Z}}}  
\newcommand{\NN}{{\cal N}}
\renewcommand{\AA}{{\cal A}}
\newcommand{\tR}{\tilde{R}}
\newcommand{\bR}{\bar{R}}
\newcommand{\nn}{\nonumber}

\newcommand{\tPhi}{{\tilde{\Phi}}}
\newcommand{\tphi}{{\tilde{\phi}}}
\newcommand{\tDe}{{\tilde{\De}}}

\def\be#1{\begin{equation} \label{#1}}
\def\bes{\end{equation}}
\def\bi{\begin{itemize}}
\def\bs{\begin{split}}
\def\es{\end{split}}

\def\ba{\begin{align}}
\def\bas{\begin{align*}}
\def\ea{\end{align}}
\def\eas{\end{align*}}

\def\R{{\hbox{\bf R}}}

\newcommand{\bea}{\begin{eqnarray}}
\newcommand{\eea}{\end{eqnarray}}
\newcommand{\beaa}{\begin{eqnarray*}}
\newcommand{\eeaa}{\end{eqnarray*}}

\def\diag{{\hbox{diag}}}

\def\nab{\nabla}
\def\nabo{\overline{\nabla}}

\def\varep{\varepsilon}

\newenvironment{proof}{\noindent {\bf Proof} }{\endprf\par}
\def \endprf{\hfill  {\vrule height6pt width6pt depth0pt}\medskip}
\def\emph#1{{\it #1}}
\def\textbf#1{{\bf #1}}

\theoremstyle{plain}
  \newtheorem{theorem}[subsection]{Theorem}
  
  \newtheorem{proposition}[subsection]{Proposition}
  \newtheorem{lemma}[subsection]{Lemma}

\theoremstyle{remark}
  \newtheorem{remark}[subsection]{Remark}

\theoremstyle{definition}
  \newtheorem{definition}[subsection]{Definition}

\include{psfig}

\begin{document}

\title[]{On the global regularity of Wave Maps
in the critical Sobolev norm}
\author{Sergiu Klainerman}
\address{Department of Mathematics, Princeton University, Princeton NJ 08544}
\email{ seri@@math.princeton.edu}

\author{Igor Rodnianski}
\address{Department of Mathematics, Princeton University, Princeton NJ 08544}
\email{ irod@@math.princeton.edu}

\vspace{-0.3in}
\begin{abstract}
We extend the recent result of T.Tao \cite{T}
to wave maps defined from the Minkowski space
$\R^{n+1}$, $n\ge 5$,  to  a  target manifold $\NN$
 which possesses a ``bounded parallelizable'' structure.
This is the case of Lie groups, homogeneous spaces
as well as the hyperbolic spaces ${\bf H}^N$.  General compact
Riemannian manifolds can be imbedded as totally geodesic
submanifolds  in bounded parallelizable
manifolds, see \cite{CZ},  and therefore are also covered, in principle,
 by   our result. Compactness of the target manifold, which 
 seemed to play an important role in \cite{T}, turns out  however to play no 
role in our discussion.
 Our proof
follows closely that of \cite{T} and is based, in particular,
on its remarkable  microlocal  gauge renormalization idea.

\end{abstract}

\maketitle

\section{Introduction}

\qquad Let $\phi: \rr^{n+1}\longrightarrow (\NN,h)$ with $(\NN,h)$
a Riemannian manifold of dimension $N$. Here $\R^{n+1}$ denotes the 
standard Minkowski space endowed with
 the metric $m=\diag(-1,1,\ldots, 1).$ We denote by $\nab$ the 
Levi-Civitta  connection  on $TN$, the tangent
bundle of $\NN$, and by $\nabo$ the  induced 
connection on  $\phi^*(T\NN)$. Recall that  the pull-back bundle 
$\phi^*(T\NN)=\cup_{x\in \rr^{n+1}}\,\,\,\{x\}\times T_{\phi(x)}\NN $ 
 is   a  vector bundle over  the  Minkowski space $\rr^{n+1}$.
The induced metric on  $\phi^{*}(T\NN) $ is defined by
$$<V,W>=V^a(\phi(x))V^b(\phi(x))<e_a,e_b>_h$$
where $V=V^a(\phi(x))e_a$, $W=W^b(\phi(x))e_b$
are sections of $\phi^*(T\NN)$ and $e_a$ is a
frame of vectorfields on $\NN$. The induced connection $\nabo$
 is defined according to the rule
$$\nabo_XV=\nab_{\phi_*X} V$$
where $X\in T(\R^{n+1})$ and $V\in \phi^*(T\NN)$.
 A map $\,\,\phi: \rr^{n+1}\longrightarrow (\NN,h)$ is said to be a wave map if 
$$ m^{\a\b}\nabo_{\pr_\b}\phi_{*}(\pr_\a)=0.$$
In local coordinates $y^I, I=1,\ldots,N$ on $\NN$ the wave maps equation takes
the familiar form
\be{1.4}
\partial^\a\pr_\a\phi^I+\Ga^I_{JK}\pr_\b\phi^J\pr_\ga\phi^K m^{\b\ga}=0.
\end{equation}
where $\Ga^I_{JK}$ are the Christoffel coefficients
of the Levi-Civitta connection $\nabla$ on $\NN$ and $\phi^I, I=1,\ldots,N$
the components of the map $\phi$ in local coordinates on $\NN$.

Let  $e_a=e_a^I\frac{\pr}{\pr y^I}$ be an  orthonormal frame of vectorfields
 and 
$\om^a=\om^a_Idy^I$ be the corresponding dual basis of 1-forms
$\om^a(e_b)=\de^a_b$. 
 Since $h(e_a,e_b)=\de_{ab}$ we infer that $h_{IJ}=\sum_a \om^a_I \om^a_J$.
Define,
\be{PHI}
\phi^a_\a=\om^a_I\pr_\a\phi^I
\end{equation}
where $\phi^I$ are the components of the map
$\phi$ relative to the local coordinates $y^I$ on $\NN$.
Clearly,
$\pr_\a\phi^I=e_a^I\phi^a_\a.
$
Given a function  $F$ on $\NN$ we  write
$\pr_\a F(\phi)=\frac{\pr}{\pr y^I} F(\phi)\pr_\a\phi^I=
\frac{\pr}{\pr y^I} F(\phi)e^I_d\phi^d_\a=F_{,d}(\phi)\phi^d_\a $
where $F_{,d}=e_d(F)$.

\qquad We easily  check that the functions $\phi^a_\a=<\pr_\a\phi,e_a>$ associated
to a wave map $\phi$ verify( see 
\cite{CZ}) the following divergence-curl system,
\bea
\pr_\b\phi^a_\a-\pr_\a\phi^a_\b&=&C^a_{bc}\,\,\phi^b_\a\,\phi^c_\b\label{1.2}\\
\pr^\a\phi^a_\a&=&-\Ga^a_{bc}\,\phi^b_\b\,\phi^c_\ga\,m^{\b\ga}\label{1.3}
\eea
where, $C^a_{bc}$ and $\Ga^a_{bc}$ are respectively  the structure and connection
coefficients of the frame,
\beaa
[e_b,e_c]&=&C^a_{bc}e_a\\
\nab_{e_b}e_c&=&\Ga^a_{bc}e_a
\eeaa
 In view of  the formula $\,\,[e_b,e_c]=\nab_{e_b}e_c-\nab_{e_c}e_b$
we infer that,
$$C^a_{bc}=\Ga^a_{bc}-\Ga^a_{cb}.$$

Since  the  frame $e_a$ is orthonormal
we have $\Gamma_{bc}^a=<\nab_{e_b} e_c\,;\,e_a> =-< e_c\,;\nab_{e_b}\,e_a> $ and
therefore
\be{1.5}
\Ga^a_{bc}=-\Ga^c_{ba}.
\end{equation}
Also,
$$\Ga^a_{bc}=\frac 1 2\bigg(C^a_{bc}+C^b_{ac}+C^c_{ab}\bigg).$$

\medskip
\noindent
{\bf Definition:}\,\,\,{\em We say that a Riemannian  manifold 
$\NN$ has a `` bounded parallelizable'' structure if  there
exists an orthonormal frame $(e_a)_{a=1}^N$ on $\NN$  relative 
to which the structure coefficients 
$C^a_{bc}$  and their frame derivatives $C^a_{bc\,, \,d_1d_2\ldots d_k}$ are 
 uniformly bounded on $\NN$. } 
\begin{remark}
There are plenty of examples of bounded parallelizable manifolds.
To start with on  any Lie group  we can construct an orthonormal 
basis of left invariant vectorfields $\, e_a\,$ relative to which
the structure constants $C^a_{bc}$ are constant\footnote{We   refer to
these  as ``constant parallelizable''.}.
The constant negative curvature manifolds ${\bf H}^N$, $N>2$
 are bounded parallelizable. Moreover ${\bf H}^2$, i.e. the hyperbolic plane, 
is a Lie group,  see  the relevant
discussion in  section 3.1 of \cite{CZ}.
In addition any compact Riemannian manifold 
can be embedded  as a totally geodesic submanifold in a bounded parallelizable
Riemannian  manifold, see \cite{CZ}.
\end{remark}
\begin{proposition}
\label{equation} Let $\NN$ be a  Riemannian manifold
and $\phi:\R^{n+1}\longrightarrow \NN$ a wave map. The 1-forms 
  $\phi^a_\a=<\pr_\a\phi,e_a>$ verify the equations, \eqref{1.2},\eqref{1.3}
as well as the system of wave equations, 
\be{1.8}
\Box\Phi=-2R_\mu\,\cdot \,\pr^\mu\Phi +E
\end{equation}
with $\Phi=(\phi^a_\a)$, $R_\mu=(R^a_{b\mu})_{a,b=1}^N$ and 
$ R^a_{b\mu}=\Ga^a_{cb}\phi^c_\mu. $ The components of $E=(E^a_\a)$
are homogeneous polynomial of degree   three  relative to  the components
 of $\Phi=(\phi^a_\a)$
with coefficients depending only on  the structure functions 
$C^a_{bc}$ and their derivatives $C^a_{bc,d}$ with respect to the frame.
\end{proposition}
\begin{remark}
It is essential to remark that the matrices $R_\mu$ are antisymmetric i.e.
\be{1.7'}
R^a_{b\mu}=-R^b_{a\mu}
\end{equation}
This is an immediate consequence of \eqref{1.5}. This shows that
the well known  ``Helein trick'' of  antisymmetrizing the form of the wave maps 
equations in the particular case when $\NN$ is a standard sphere, trick
which plays a fundamental role in \cite{T}, is due  in fact to a general feature
of the connection coefficients on \emph{ any  Riemannian manifold}, expressed
relative to  \emph{ orthonormal frames}. 
\end{remark}
\begin{proof}:

\qquad Differentiating \eqref{1.2}
and using \eqref{1.3} we derive:
\beaa
\pr^\b\pr_\b\phi^a_\a&=&-\pr_\a(\Ga^a_{bc}\,\phi^b_\mu\,\phi^c_\nu \,m^{\mu\nu})
+\pr^\b(C^a_{bc}\,\phi^b_\a\,\phi^c_\b)\\
&=&m^{\mu\nu}\bigg(-\Ga^a_{bc}\,(\pr_\a\phi^b_\mu)\,\phi^c_\nu
-\Ga^a_{bc}\,\phi^b_\mu(\,\pr_\a\phi^c_\nu\,) 
+C^a_{bc}(\,\pr_\mu\phi^b_\a)\,\phi^c_\nu \bigg)\\
&+& C^a_{bc}\,\phi^b_\a(\,\pr^\b\phi^c_\b\,)
-\pr_\a(\Ga^a_{bc})\,\phi^b_\mu\,\phi^c_\nu
\,m^{\mu\nu}+\pr^\b(C^a_{bc})\,\phi^b_\a\,\phi^c_\b
\eeaa
Setting  $A^a_\a=-m^{\mu\nu}\bigg(\Ga^a_{bc}\,(\pr_\a\phi^b_\mu)\,\phi^c_\nu
+\Ga^a_{bc}\,\phi^b_\mu(\,\pr_\a\phi^c_\nu\,) 
-C^a_{bc}(\,\pr_\mu\phi^b_\a)\,\phi^c_\nu \bigg)$  and using 
\eqref{1.2} we write
\beaa
A^a_\a&=&-m^{\mu\nu}\bigg(\Ga^a_{bc}\,(\pr_\mu\phi^b_\a
+C^b_{mn}\phi^m_\a\phi^n_\mu)\,\phi^c_\nu+
\Ga^a_{bc}\,\phi^b_\mu(\pr_\nu\phi^c_\a+ C^c_{mn}\phi^m_\a\phi^n_\nu)\\
&-&C^a_{bc}(\,\pr_\mu\phi^b_\a)\,\phi^c_\nu \bigg)\\
&=&-m^{\mu\nu}\bigg(\Ga^a_{bc}+\Ga^a_{cb}-C^a_{bc}\bigg)\phi^c_\nu\pr_\mu\phi^b_\a
+m^{\mu\nu}\Ga^a_{bc}\bigg(C^b_{mn}\phi^m_\a\phi^n_\mu\phi^c_\nu+
C^c_{mn}\phi^m_\a\phi^n_\nu\phi^b_\mu\bigg)
\eeaa
or since
$\Ga^a_{bc}+\Ga^a_{cb}-C^a_{bc}=\Ga^a_{bc}+\Ga^a_{cb}-(\Ga^a_{bc}-\Ga^a_{cb})=
2\Ga^a_{cb}, $
$$
A^a_\a=-2m^{\mu\nu}\Ga^a_{cb}\phi^c_\nu\pr_\mu\phi^b_\a
+m^{\mu\nu}\Ga^a_{bc}\bigg(C^b_{mn}\phi^m_\a\phi^n_\mu\phi^c_\nu+
C^c_{mn}\phi^m_\a\phi^n_\nu\phi^b_\mu\bigg)
$$
Therefore,
\beaa
\pr^\b\pr_\b\phi^a_\a&=&-2m^{\mu\nu}\Ga^a_{cb}\phi^c_\nu\pr_\mu\phi^b_\a+m^{\mu\nu}\Ga^a_{bc}\bigg(C^b_{mn}\phi^m_\a\phi^n_\mu\phi^c_\nu+
C^c_{mn}\phi^m_\a\phi^n_\nu\phi^b_\mu\bigg)\\
&-&C^a_{bc}\,\phi^b_\a(\,\Ga^c_{mn}\phi^m_\mu\phi^n_\nu m^{\mu\nu}\,)
-\Ga^a_{bc,d}\,\phi^b_\mu\,\phi^c_\nu\phi^d_\a
\,m^{\mu\nu}+C^a_{bc,d}\,\phi^b_\a\,\phi^c_\b\phi^d_\b
\eeaa
Finally we write
$$
\Box\phi^a_\a=-2\Ga^a_{cb}\,\phi^c_\nu\,\pr_\mu\phi^b_\a m^{\mu\nu} +E^a_\a
$$
where $E=E^a_\a$ denote  the terms above  which are cubic in $\Phi=(\phi^a_\a)$
with coefficients depending only on  the structure functions
$C^a_{bc}$ and their derivatives $C^a_{bc,d}$ with respect to the frame.
This is precisely the equation \eqref{1.8}.
\end{proof}
\qquad We study the evolution of  wave maps subject to the initial value problem
\be{in1}
\phi(0)=\varphi, \qquad \pr_t\phi(0)=\psi=\psi_0^ae_a
\end{equation}
$\varphi$ is an arbitrary smooth  map  defined from $\R^n$ 
with values in $\NN$ and $\psi=\psi_0^ae_a$ and arbitrary  smooth map
 from $\R^{n+1}$ to $TN$.  Let $\varphi^a_i=<\pr_i\varphi,e_a>$. 
\begin{definition}\,\,
We shall say that the initial data $\phi[0]=(\varphi,\psi)$
belongs to the Sobolev space  ${\dot H}^s(\R^{n})$, resp. $H^s(\R^{n})$,  if all
components 
$\varphi^a_i,\,\, \psi^a_i$ belong to the space ${\dot H}^{s-1}(\R^{n})$,
resp. $H^{s-1}(\R^{n})$. We write
$$\|\phi[0]\|_{\dot{H}^s}=\sum_{a,i}\bigg(\|\varphi^a_i\|_{\dot{H}^{s-1}}+
\|\psi^a_i\|_{\dot{H}^{s-1}}\bigg)
$$ and similarly for $\|\phi[0]\|_{{H}^s}$.
\label{Hs}
\end{definition}
 
 We are now ready to 
state our main theorem.

{\bf Main Theorem}\quad {\em  Let $ \NN$ be a Riemannian manifold
endowed with  a bounded  parallelizable structure. Assume
$n\ge 5$ and that the initial data $\phi[0]=(\varphi, \psi=\psi^a_i e_a)$ is
 in $H^s$ for some $\,\frac{n}{2}<s$. We make also the
 critical smallness  assumption:
$$\|\phi[0]\|_{\dot{H}^{\frac{n}{2}}}\le \varep$$
Then the wave map $\phi$ with initial data $\phi[0]$ can be 
uniquely continued in  $H^s$ norm
globally in time.
}

\qquad The proof of the Main Theorem relies on a local well-posedness
result in $H^s$, $s>\frac{n}{2}$. We state the  precise result below:
\begin{theorem} Assume that the initial data $\phi[0]\in H^{s}(\R^n)$
for some  $s\ge s_0>\frac{n}{2}$. There exists a
$T>0$, depending only on the size of $\|\phi[0]\|_{H^{s_0}}$,  and a unique
solution
$\phi$ of the system \eqref{1.2} ,\,\eqref{1.3} defined on  the slab
 $[0,T]\times \R^n$ verifying, 
$$\|\phi[t]\|_{H^s}\le  C \|\phi[0] \|_{H^s}$$
for all $t\in [0,T]$ and $C$ a constant depending only on $T$ and $s$.
\label{local}
\end{theorem}
\bigskip

\qquad Strictly speaking such a sharp local existence result
for  div-curl systems of type \eqref{1.2}, \eqref{1.3}
does not exist in the literature. Nevertheless we are confident
that the methods discussed in \cite{KSe} in connection to 
a special model problem related to \eqref{1.2},\eqref{1.3}( see also
\cite{KMa}) do apply\footnote{In \cite{KMa} one proves
a related result for a model problem in dimension $n=3$.
The proof was vastly simplified  and extended to all dimensions $n\ge 3$
in \cite{Se} and \cite{KSe}. The higher dimensional case is in fact a lot simpler. 
}. A proof of this fact will appear elsewhere. Alternatively we can avoid
Theorem \ref{local} and rely instead on the known sharp local existence result
for the Wave Maps system  written in local coordinates \eqref{1.4}, see
\cite{KSe}
 and the references therein. Indeed any  $H^s$ data\footnote{According to
the  global  definition \ref{Hs}.},
$s>\frac{n}{2}$,  is also $H^s$ with respect
to local coordinates on
$\NN.$ Using the finite propagation speed property of wave equations we can
therefore construct a local in time  $H^s$ solution for the system
\eqref{1.2}-\eqref{1.3} which is unique, as a solution of \eqref{1.4},  
 in any local chart on $\NN$.  This is the solution for  which the
Main Theorem applies.

\qquad For simplicity
we shall present the proof of the Main Theorem in the particular case
of constant parallelizable target  manifolds. The general case
complicates matters only in so far as the number of terms we need to treat
is larger, there are however no conceptual differences. We shall thus
assume that the $C^a_{bc}$ are constant and  indicate
whenever needed what additional steps  are required to  treat the general case.
 
\emph{Aknowledgement:}\,\,\, We would like to thank D. Christodoulou and 
T. Tao for valuable comments.
 
\section{Notation, Strichartz estimates and main proposition}
\qquad   We use the Littlewood-Paley notation of \cite{T}.
Thus, for a function $\phi(t,x)$ we denote the projections
$P_{\le k}\phi(t,x)= \int e^{ix\cdot \xi}\chi(2^{-k}\xi)\phi\,\hat{}(t,\xi) d\xi$ 
where $\phi\,\hat{}(t,\xi)$ is the space Fourier transform of $\phi$ and 
$\chi(\xi)=\eta(\xi)-\eta(2\xi)$  with $\eta$  a non-negative smooth  bump function
supported on $|\xi|\le 2$  and equal to $1$ on the ball $|\xi|\le 1$. Therefore
$\chi(\xi)$ is supported in $\{\frac{1}{2}\le |\xi|\le 2\}$ and $\sum_{k\in {\bf
Z}}\chi(2^{-k}\xi)=1$ for all $\xi\ne 0$.  We also define $P_k=P_{\le k}-P_{<k}$.
Also for any interval
$I\subset {\bf Z}$ we define  $P_I$ in an obvious fashion, see \cite{T}.

  \qquad Following \cite{T} we introduce the  notation
\be{1.9}
\|\Phi\|_{S_k}=\sup_{q,r\in \AA}2^{k(\frac{1}{q}+\frac{n}{r}-1)}
\bigg(\|\Phi\|_{L_t^qL_x^r}+2^{-k}\|\pr_t\Phi\|_{L_t^qL_x^r}\bigg)
\end{equation}
where $\AA=\{(q,r)/ 2\le q,r\le \infty, \,\,\frac{1}{q}+\frac{n-1}{2r}\le
\frac{n-1}{4}\}$ is the set of admissible Strichartz exponents.
Recall that,
\begin{theorem}
\label{Strich}
For any fixed integer $k$ and $\phi(t,x)$ a function on $\R\times \R^n$
such that the support of $\hat{\phi} \,(t,\xi)$ is included in the 
dyadic region $2^{k-1}\le |\xi|\le 2^{k+1}$ we have the estimate,
$$\|\phi\|_{S_k}\lesssim 
\|\phi(0)\|_{H^{\frac{n-2}{2}}}+
\|\pr_t\phi(0)\|_{H^{\frac{n-4}{2}}}+2^{k\frac{n-4}{2}}
\|\Box \phi\|_{L_t^1L_x^2}.$$
\end{theorem}
In what follows we recall the definition
of frequency envelope given in \cite{T}.
\begin{definition} A frequency envelope 
is an $l^2$  sequence $c=(c_k)_{k\in {\bf Z}}$ verifying
\be{env}
c_{k}\lesssim 2^{\si|k-k'|} c_{k'}
\end{equation}
for all $k,k'\in {\bf Z}$. Here $\si$ is a fixed
positive constant; as in \cite{T} we take $0<\si<\frac{1}{2}$.
In addition we shall also need $0<\si<\frac{n-4}{4}$.

\qquad We say that the  $\dot{H}^{s}$   norm of a function $f $ on $\R^n$ lies
underneath an envelope $c$ if, for all $k\in {\bf Z}$, 
$$\|P_kf \|_{\dot{H}^{s}}\le c_k.$$
We shall write $ f<<_s c$ or simply $f<<c$ when there is no danger 
of confusion.
Recall, see \cite{T} section 3, that if $\|f\|_{\dot{H}^s}\le \varep$
then there exists an envelope $c\in l^2$ such that
 $\|c\|_{l^2}\lesssim \varep$ and $f<<_s c$. Indeed we can  simply take,
$c_k=\sum_{k'\in {\bf Z}}2^{-\si|k-k'|}\|P_k\|_{\dot{H}^s}$.
\end{definition}
\begin{definition} Fix $0<\si< \min(\frac{1}{2}, \frac{n-4}{4})$  and $c$ a 
frequency envelope. We say that the initial data
$\phi[0]=\bigg(\phi(0)=\varphi, \pr_t\phi(0) =\psi=\psi^a_0 e_a\bigg)$ lies
underneath
$c$ if, relative to our frame $e_a$  we have  for all $k\in {\bf Z}$,
$$\| P_k\phi[0]\|_{\dot{H}^{\frac{n}{2}}}\le c_k.$$
We shall use the short hand notation  $\phi[0]<<c.$
\label{defin}
\end{definition}
\qquad   Following the same arguments as in section 3 of \cite{T}
we can reduce the proof of our  main theorem to the following\footnote{Our
main proposition below corresponds to Proposition 3.3 in \cite{T}. The
reduction relies on Theorem \ref{local}. }:

\begin{proposition}( Main Proposition)\quad Let $c$ be a frequency
envelope\footnote{verifying \eqref{env} with $\si< \min(\frac 1 2,
\frac{n-4}{4})$.} with
$\|c\|_{l^2}\le
\varep$,
$0<T<\infty$ and
$\Phi=(\phi^a_\a)$ verify the equations
 \eqref{1.2}, \eqref{1.3},  and therefore also \eqref{1.8}. Assume
  that,  according to definition  \ref{defin}, the initial data verifies
 the smallness condition $\phi[0]<<c$.  Assume also the bootstrap assumption,
\be{1.10}
\|P_k\Phi\|_{S_k([0,T]\times\R^n)}\le 2C c_k
\end{equation}
for all $k\in \Z$.
Then in fact, for sufficiently small $\varep$, and all $k\in \Z$,
\be{1.11}
\|P_k\Phi\|_{S_k([0,T]\times\R^n)}\le C c_k.
\end{equation}
\end{proposition}
Returning to the definition \eqref{1.9} we make explicit all   the 
 useful estimates contained in the  bootstrap assumption \eqref{1.10},
\beaa
\|P_k\Phi\|_{L_t^2L_x^{\frac{2(n-1)}{n-3}}}+2^{-k}\|\pr_t
 P_k\Phi\|_{L_t^2L_x^{\frac{2(n-1)}{n-3}}}&\le&2^{\frac k 2
-\frac{nk}{2}+\frac{nk}{n-1}}\cdot (2Cc_k)\\
\|P_k\Phi\|_{L_t^2L_x^4}+
2^{-k}\|\pr_t P_k\Phi\|_{L_t^2L_x^4}&\le& 2^{\frac k 2-\frac{nk}{4}}\cdot(2Cc_k)\\
\|P_k\Phi\|_{L_t^2L_x^{n-1}}+2^{-k}\|\pr_t P_k\Phi\|_{L_t^2L_x^{n-1}}&\le&
2^{\frac k 2-\frac{nk}{n-1}}\cdot(2Cc_k)\\
\|P_k\Phi\|_{L_t^2L_x^{\infty}}+2^{-k}\|\pr_t P_k\Phi\|_{L_t^2L_x^{\infty}}&\le&
2^{\frac k 2}\cdot(2Cc_k)\\
\|P_k\Phi\|_{L_t^\infty L_x^{2}}+2^{-k}\|\pr_t
P_k\Phi\|_{L_t^\infty L_x^{2}}&\le& 2^{\frac{k}{2}(2-n)}\cdot(2Cc_k)\\
\|P_k\Phi\|_{L_t^\infty L_x^{\infty}}+2^{-k}\|\pr_t
P_k\Phi\|_{L_t^\infty L_x^{\infty}}&\le& 2^k\cdot(2Cc_k)\\
\|P_k\Phi\|_{L_t^4 L_x^{2(n-1)}}+2^{-k}\|\pr_t P_k\Phi\|_{L_t^4L_x^{2(n-1)}}&\le&
2^{\frac{3k}{4}-\frac{nk}{2(n-1)}}\cdot(2Cc_k)\\
\eeaa
\begin{lemma} The assumptions \eqref{1.10} imply 
\be{le}
\|\square P_k\Phi\|_{S_k}\lesssim 2^{2k} C c_k
\end{equation}
\label{lbox}
\end{lemma}
\begin{remark} The Lemma \ref{lbox} seems morally right yet somewhat
involved to prove it in details.
In the applications below, see also \cite{T},  we shall   only need the
estimate
\be{box}\|\square P_k\Phi\|_{L_t^2L_x^{n-1}}\lesssim
 2^{k( 2+\frac 1 2-\frac{n}{n-1})} C c_k.
\end{equation}
We prove this estimate  at the end of section 3.
\end{remark}

 \qquad In view of the scale invariance of both  our equations and the smallness 
condition $\phi[0]<<c$ 
it suffices  to prove \eqref{1.11}\,  for $\, k=0$. Let  $\Psi=P_0\Phi$.
We need to prove that,
\be{1.12}
\|\Psi\|_{S_0([0,T]\times\R^n)}\le C c_0
\end{equation}

\qquad To  prove \eqref{1.12}  we would like  to apply Theorem \ref{Strich}
to the equation  obtained by applying the projection $P_0$.
to  \eqref{1.8} i.e.,
$$\square \Psi=P_0(R_\mu \cdot \pr^\mu \Phi+E). $$
A straightforward application of the Strichartz 
inequalities will not work however. Indeed according
to Theorem \ref{Strich}
$$\|\Psi\|_{S_0}\le c_0 +\|P_0(R_\mu \cdot \pr^\mu \Phi+E)\|_{L_t^1L_x^2}$$
\qquad The cubic term $E$ presents no difficulty, the problem comes up
when we try to estimate $ P_0(R_\mu \cdot \pr^\mu \Phi)$ more precisely
the part of it which corresponds to the interaction between low frequencies
of $R$ and frequencies of $\Phi$  comparable to  those of $\Psi$. More
precisely the most dangerous terms are of the form $\tR \cdot \pr \Psi$
with $\tR=P_{\le -10}R$. To estimate $\|\tR \cdot \pr\Psi\|_{L_t^1L_x^2}$
relative to the available Strichartz norms we are forced to take 
$\Psi$ in the energy  norm $L_t^\infty L_x^2$.  This leaves  us  with
$\tR $ in the norm $L_t^1L_x^\infty$ for which we don't have any
Strichartz estimates. It is  precisely this  difficulty  which  led Tao to
introduce his remarkable renormalization idea which we reproduce below in section
5. 

\qquad The organization of the  paper follows closely  that of \cite{T}.
In the next section we reduce the proof of the main proposition
to estimates for the  linearized equation:
$$ \square \Psi=-2\tR_\mu \cdot \pr^\mu \Psi.$$ This corresponds to 
isolating the worst part of 
$P_0(R_\a\cdot\pr^\a\Phi)$ to which we have alluded above. 
 In section 4, which represents the main
  contribution of this
paper, we show how to replace the term $\tR_\mu$ by the perfect
derivative $\pr_\mu \tDe$ of an antisymmetric  potential $\tDe$.  This fact
plays a crucial role in carrying out Tao's renormalization procedure
 in section 5.

 \qquad We shall use, throughout the paper,  Tao's  convention to call an
acceptable error any function, or matrix valued function, $F$ on $[0,T]\times \R^n$
such that 
\be{error}
\|F\|_{L_t^1L_x^2( [0,T]\times \R^n)}\le C^3\varep c_0
\end{equation}

\section{Reduction to a linear equation} 
\begin{proposition} The matrix valued function $P_0\Phi=\Psi$
verifies the equation
\be{1.13}
\square \Psi=-2\tR_\mu \cdot \pr^\mu \Psi + \mbox{error} 
\end{equation}
where $\tR_\mu=P_{\le -10} R_\mu=\Ga \cdot \tPhi_\mu $ and 
$\tPhi_\a=P_{\le -10}\Phi_\a$. Here $\mbox{``error''} $
refers to  an acceptable error term in the sense of \eqref{error}.
\label{Plin}
\end{proposition}
\begin{remark} Written in components $\Psi=(\psi^a_\a)$ with
$\psi^a_\a=P_0\phi^a_\a$ the equation  \eqref{1.13}
has the form 
$$\square \psi^a_\a=-2\tR^a_{b\nu} \cdot \,\pr_\mu\psi^b_\a
m^{\mu\nu}+\mbox{error},$$
 where $\tR^a_{b\nu}=\Ga^a_{bc}\tphi^c_\nu $ and   
$\tphi^c_\nu=P_{\le -10}\phi^c_\nu$. Observe that the 
$N\times N $ matrices $\tR_\mu$ are antisymmetric i.e. $\tR_\mu^{\,t}=-\tR_\mu$.
\end{remark}
\begin{proof}:  We start with the equation \eqref{1.8} to which 
we apply the projection $P_0$. Therefore,
$$\square\Psi= P_0(R_\mu \cdot \pr^\mu \Phi +E)$$
The proof of Proposition \ref{Plin} is an immediate consequence of the following
Lemmas
\begin{lemma} We have,
$$P_0(R_\mu \cdot \pr^\mu \Phi)=\tR_\mu\cdot \pr^\mu\Psi+
\mbox{error}$$
where $\tR_\mu=P_{\le -10} R_{\mu}$.
\label{le3.1}
\end{lemma}
\begin{lemma} The term $P_0E$ is an  acceptable error term.
\label{le3.2}
\end{lemma}
\end{proof}
\qquad We  sketch below the proof of Lemma  \ref{le3.1}. Lemma
\ref{le3.2} is easier and can be proved in a similar manner\footnote{The proof
of the Lemma  requires only the   boundedness of the structure coefficients
 $C$ and their frame derivatives. Cubic
terms in $\Phi$ are easy to treat by the available
Strichartz inequalities.}.
We start by decomposing $R_\mu=\sum_k P_k R_\mu=\sum_k R_{\mu,k}$ and 
$\Phi=\sum_k P_k\Phi=\sum \Phi_k$.  Thus,
$$P_0(R_\mu \cdot \pr^\mu \Phi)=P_0(E_1+E_2+E_3 +E_4+E_5+E_6)$$
where,
\beaa
E_1&=& \sum_{\max(k_1,k_2)>10\,,\, |k_1-k_2|\le 5}R_{\mu,k_1}\cdot\pr^\mu
\Phi_{k_2}\\
E_2&=& \sum_{\max(k_1,k_2)>10\,,\, |k_1-k_2|> 5}R_{\mu,k_1}\cdot\pr^\mu
\Phi_{k_2}\\
E_3&=&  (P_{\le -10}R_{\mu})\cdot \pr^\mu
P_{\le -10}\Phi\\
E_4&=&  (P_{(-10,10)}R_{\mu})\cdot \pr^\mu
P_{(-10,10)}\Phi\\
E_5&=&  \tR_{\mu}\cdot\pr^\mu
P_{(-10,10)}\Phi\\
E_6&=&  (P_{(-10,10)}R_{\mu})\pr^\mu
P_{<-10}\Phi
\eeaa
 Recall that  the matrices $R_\mu$ are products
between the constant matrices\footnote{In the general case
of a bounded parallelizable manifold one has to take into account
the additional commutator terms. The commutators generate   additional
powers of $\Phi$  and therefore can be treated as  easy error terms.} $\Ga$ and
$\Phi$. Thus each 
$P_k R_\mu$ can be estimated in the same way as $P_k\Phi=\Phi_k$ according to
\eqref{1.10}. 
  Using \eqref{1.10} and the
envelope property
\eqref{env},
\beaa
\|E_1\|_{L_t^1L_x^2}&\le& \sum_{\max(k_1,k_2)>10\,,\, |k_1-k_2|\le 5}
\|R_{\mu,k_1}\|_{L_t^2L_x^4}\cdot\|\pr^\mu\Phi_{k_2}\|_{L_t^2L_x^4}\\
&\lesssim& C^2\sum_{k\ge 5} 2^k 2^{k-\frac{nk}{2}} c_k^2\lesssim
C^2c_0^2\sum_{k\ge 5}
2^k \,2^{k-\frac{nk}{2}}\, 2^{2\si k}\\
&\lesssim& C^2\, c_0^2\,\sum_{k\ge 5}2^{k(2+2\si-\frac{n}{2})}\lesssim C^2\, c_0^2
\eeaa
provided that $\si<\frac{n-4}{4}$.

\qquad Clearly  $P_0 E_2=0, P_0E_3= 0$.  The term $E_4$ is easy
to estimate; it contains  only a finite number of terms,
$$
\|E_4\|_{L_t^1L_x^2}\le \|P_{(-10,10)} R\|_{L_t^2L_x^4}
\|\pr P_{(-10,10)}\Phi\|_{L_t^2L_x^4}\lesssim C^2 c_0^2.
$$
\qquad For $E_6$ we write,
\beaa
\|E_6\|_{L_t^1 L_x^2}&\lesssim & \|\pr P_{\le
-10}\Phi\|_{L_t^2L_x^{\frac{2(n-1)}{n-3}}}
\|\pr P_{(-10,10)}\Phi\|_{L_t^2L_x^{n-1}}\\
&\lesssim&C^2c_0\sum_{k\le -10}2^{k(\frac 3 2+ \frac{1}{n-1})}c_k\\
&\lesssim&C^2c_0^2\sum_{k\le -10}2^{k(\frac 3 2+ \frac{1}{n-1}-\si)}\lesssim
C^2c_0^2
\eeaa
\qquad It remains to consider the term $P_0E_5=P_0\bigg(\tR_{\mu}\cdot\pr^\mu
P_{(-10,10)}\Phi\bigg)$.
 We use the standard commutator
inequality for functions $f,g$ in $\R^n$, see Lemma 4.3 in \cite{T},
$$\|P_0(fg)-fP_0g\|_{L^r}\lesssim \|\nabla f\|_{L_p}\|g\|_{L^r}
$$
Taking $r=2$, $p=n-1$ and $q=\frac{2(n-1)}{n-3}$  and proceeding
precisely as for $E_6$ we derive,
\beaa
\|P_0\big(\tR_{\mu}\cdot\pr^\mu
P_{(-10,10)}\Phi\big)-\tR_{\mu}\cdot\pr^\mu\Psi\|_{L_t^1L_x^2}&\lesssim&
\|\pr \tR_\mu\|_{L_t^2L_x^{n-1}}\|\pr
P_{(-10,10)}\Phi\|_{L_t^2L_x^{\frac{2(n-1)}{n-3}}}\\
&\lesssim&C^2c_0^2\sum_{k\le -10}2^{k(\frac 3 2+ \frac{1}{n-1}-\si)}\lesssim
C^2c_0^2
\eeaa
as desired.
\medskip

\qquad In the remaining part of this section we shall
 sketch the proof of  Lemma \ref{lbox}. More precisely 
we derive the estimate \eqref{box}, all other estimates
can be derived in a similar manner. By scale invariance
it suffices to prove \eqref{box} for $k=0$. In other words
we have to prove,
\be{boxx}
\|\square \Psi\|_{L_t^2L_x^{n-1}}\lesssim C c_0
\end{equation}
According to Lemma \ref{le3.1} it suffices to prove that
$\|\tR_\mu\cdot \pr^\mu\Psi\|_{L_t^2L_x^{n-1}}\lesssim C c_0.$
 Using \eqref{1.10} we derive ,
\beaa
\|\tR_\mu\cdot \pr^\mu\Psi\|_{L_t^2L_x^{n-1}}&\le &
\|\tR\|_{L_t^4 L_x^{2(n-1)}} \|\pr \Psi\|_{L_t^4 L_x^{2(n-1)} }\\
&\lesssim& 
C^2c_0\sum_{k\le -10}2^{k(\frac{3}{4}-\frac{n}{2(n-1)})}c_k\\
&\lesssim& 
C^2c^2_0\sum_{k\le -10}2^{k(\frac{3}{4}-\frac{n}{2(n-1)}-\si)}\lesssim C^2c_0^2
\eeaa 
as desired.
\section{ Can Replace $\tR_\mu$ by $\pr_\mu \tDe$  } 
\qquad This   reduction step is   the main contribution of our  paper.
In order to apply Tao's renormalization  procedure we express 
$\tR_\mu$ in terms of the space-time gradient  of a  potential $\tDe$
plus terms which lead to error terms. More precisely,
\begin{proposition} The matrix valued function  $\Psi$ verifies
an equation of the form,
\be{1.14}
\Box\Psi=-2\pr_\mu \tDe\cdot \pr^\mu\Psi+\mbox{error}
\end{equation}
where the potential $\tDe$ verifies the following 
properties:

\medskip
\noindent
i.) \,\,The   $N\times N$ matrix $\tDe$ is antisymmetric i.e.  $\tDe^t=-\tDe$.
The  space Fourier transform of each component of $\tDe$ is supported
in $|\xi|\le 2^{-10}$.

\medskip
\noindent
ii.)\,\, The following estimates hold for any $\tDe_k=P_k\tDe$:
\bea
\|\tDe_k\|_{S_k}&\lesssim& 2^{-k}Cc_k \label{1.15}\\
\|\partial \tDe_k\|_{S_k}&\lesssim& Cc_k \label{1.16}
\eea
Also,
\be{1.17}
\|\Box  \tDe_k\|_{S_k}\lesssim 2^kCc_k
\end{equation}

\medskip
\noindent
iii.) \,\, Set $\bR_\mu=\tR_\mu-\pr_\mu\tDe$.  The following estimates hold
for  all $P_k\bR$,

\bea
\|P_k\bR\|_{L_t^1L_x^\infty}&\lesssim&C^2 c_k^2\label{1.18}\\
\|P_k\bR\|_{L_t^\infty L_x^\infty}&\lesssim& 2^k C^2c_k^2\label{1.19}
\eea
\label{step2}
\end{proposition}

\begin{proof}:\,\,
We start with the equation 
$$\square\Psi=-2\tR_\mu\pr^\mu\Psi+\mbox{error}
$$ 
We need to find the potential $\tDe$ such that
$\bR_{\mu}=\tR_{\mu}-\pr_\mu\tDe\,\,$ is small.
Clearly\footnote{ In the general case we have 
additional terms of the form $\pr \Ga(\phi)\cdot \Phi$. 
These have  the form $C'\cdot \Phi\cdot\Phi$ with $C'$
the first frame derivatives of the structure coefficients. They
are therefore similar to the terms $M\cdot\Phi\cdot\Phi$ we treat in the text.}
$$
\pr_{\nu} \bR_{\mu}-\pr_{\mu} \bR_{\nu}=
\pr_{\nu} \tR_{\mu}-\pr_{\mu} \tR_{\nu}=P_{\le-10}\big(\Ga\cdot(\pr_{\nu}
\phi_{\mu}-\pr_{\mu}
\phi_{\nu})\big)
$$
Thus according to the equation \eqref{1.2} and the constancy of 
the structure and connection coefficients,
$$\pr_{\nu} \bR_{\mu}-\pr_{\mu} \bR_{\nu}=P_{\le -10}\big(  M\cdot
\Phi\cdot\Phi\big).$$
 with $M \approx C^2$  a matrix whose entries are  quadratic in $C_{bc}^a$.
Henceforth\footnote{In the general case of a bounded parallelizable manifold one
would have an  additional commutator term which contributes, roughly speaking,  a
cubic term in $\Phi$.  },
 \beaa
\pr_{\nu}( P_k\bR_{\mu})-\pr_{\mu}(P_k \bR_{\nu})&=&M\cdot P_k(\Phi\cdot\Phi)
=E_1+E_2\\
E_1&\approx& M\cdot \sum_{k'< k-1}\Phi_{k'}\cdot\Phi_k\\
E_2&\approx& M\cdot \sum_{k_1,k_2\ge  k,\,\, |k_1-k_2|\le
2}\Phi_{k_1}\cdot\Phi_{k_2}
\eeaa
We now estimate, with the help of \eqref{1.10} and \eqref{env}
with $\si< \frac 1 2$.
\beaa
\|E_1\|_{L_t^1L_x^\infty}&\lesssim &\|\Phi_k\|_{L_t^2 L_x^\infty}\sum_{k'<
k}\|\Phi_{k'}\|_{L_t^2 L_x^\infty}\\
&\lesssim &2^{k/2} Cc_k\sum_{k'<k}2^{k'/2} C c_{k'}\lesssim 2^{k/2} C^2c_k^2
\sum_{k'<k}2^{k'/2}2^{\si(k-k')}\\
&\lesssim&2^k C^2c_k^2
\eeaa
Also, proceeding in the same way,
\beaa
\|E_2\|_{L_t^1L_x^\infty}&\lesssim &\sum_{k_1\ge k}\,\,\sum_{|k_2-k_1|\le
2}\|\Phi_{k_1}\|_{L_t^2 L_x^\infty}
\|\Phi_{k_2}\|_{L_t^2 L_x^\infty}\\
&\lesssim &\sum_{k_1\ge
k}\,\,\sum_{|k_2-k_1|\le 2}2^{k_1/2}C c_{k_1}2^{k_2/2}C c_{k_2}\\
&\lesssim &2^k C^2c_k^2
\eeaa
Therefore all the components of the exterior derivative $F_{(k)}= d(P_k\bR)$ verify
the estimates
\be{2.1}
\|F_{(k)}\|_{L_t^1L_x^\infty}\lesssim 2^k C^2c_k^2
\end{equation}
Proceeding in precisely the same manner we find that
\be{2.1'}
\| F_{(k)} \|_{L_t^\infty L_x^\infty}\lesssim 2^{2k} C^2c_k^2
\end{equation}
\qquad We define $\tDe_k$ by requiring that the spatial components
of $P_k\bR$ verify the equation,
\be{2.2}
\pr^i(P_k\bR_i)=0
\end{equation}
Consider now the divergence -curl system,
\beaa
\pr_i(P_k\bR_j)-\pr_i(P_k\bR_j)&=&F_{(k)ij}\\
\pr^i(P_k\bR_i)=0
\eeaa
By standard elliptic estimates, taking into account 
the fact that the Fourier support of $P_k\bR$   is included 
in the dyadic region $|\xi|\approx 2^k$  and using 
\eqref{2.1} we infer that,
$$\|P_k\bR_i\|_{L_t^1 L_x^\infty}\lesssim 2^{-k}\|F_{(k)}\|_{L_t^1 L_x^\infty}
\lesssim C^2c_k^2$$
On the other hand we  also have good estimates for 
$F_{(k)0i}=\pr_t P_k\bR_i-\pr_i P_k\bR_0$. In view of the divergence condition
$\pr^i(P_k\bR_i)=0$ we derive $\nabla^2 P_k\bR_0=-\pr^i F_{(k)0i}$,
 with $\nabla^2$ the Laplacean  in $\R^n$, $\nabla^2=\sum_{i=1}^n\pr_i^2$.
Therefore using standard elliptic estimates and \eqref{2.1} we infer that,
$$\|P_k\bR_0\|_{L_t^1 L_x^\infty}\lesssim 2^{-k}\|F_{(k)}\|_{L_t^1 L_x^\infty}
\lesssim C^2c_k^2$$
We have thus derived the estimate \eqref{1.18}. The estimate \eqref{1.19}
follows in the same manner from \eqref{2.1'}. We now estimate 
$\tDe$. We first observe that the divergence equation
$\pr^i(P_k\bR_i)=0$ takes the form $\nabla^2\tDe_k=\pr^i(P_k\tR_i)$.
This uniquely defines $\tDe_k$ and  we have,
$$\|\tDe_k\|_{S_k}\lesssim 2^{-k}\| P_k \tR\|_{S_k}
\lesssim 2^{-k} \| P_k \Phi\|_{S_k}\lesssim 2^{-k}C c_k.$$
which gives \eqref{1.15} and \eqref{1.16}. To prove \eqref{1.17}
we write $\nabla^2\square\tDe_k=\pr_i (P_k \square\tR_i)$.
Therefore, in view of \eqref{le},\,\,
$\|\square \tDe_k\|_{S_k}\lesssim 2^{-k}\|\square P_k \tR\|_{S_k}
\lesssim 2^{-k}\|\square P_k\tilde{\Phi}\|_{S_k}\lesssim 2^kC c_k
$
establishing  the estimate \eqref{1.17}.

\qquad To end the proof of Proposition \eqref{step2} it remains 
to observe that since each $\tDe_k$
is antisymmetric so is  the  $\tDe=\sum_{k\le -10}\tDe_k$.
 We  also need to check that the terms 
$\bR_\mu\pr^\mu \Psi$ generated when we pass from
the equation \eqref{1.13} to \eqref{1.14} are indeed error terms.
We have, using \eqref{1.18}
\beaa
\|\bR_\mu\pr^\mu \Psi\|_{L_t^1L_x^2}&\le &
\|\bR\|_{L_t^1L_x^\infty}\|\Psi\|_{L_t^\infty L_x^2}\lesssim
Cc_0 \|\bR\|_{L_t^1L_x^\infty}\\
&\lesssim& c_0 C\sum_{k\le -10}\|\bR_k\|_{L_t^1L_x^\infty}\le c_0 C 
C^2\sum_{k\le -10}c_k^2\\
&\lesssim& c_0C^3 \varep
\eeaa
as desired.

\end{proof}

\section{Tao's  renormalization procedure}
\qquad This  last   step in our proof is a straightforward implementation
of Tao's  renormalization procedure. We  repeat below  the main arguments
 in his construction.

\qquad Let $M$ be a large integer, depending on $T$, which will be chosen
below. Define the real $N\times N$ matrix  valued  function  $U$
to be 
$$U=I+\sum_{-M<k\le -10}U_k\\$$
with  the $U_k$ defined  inductively as follows,
\bea
U_k&=&0\qquad \mbox{for all\,\,\,\,} k<-M\nonumber\\
U_M&=&I\nonumber \\
U_k&=&\tDe_k\cdot U_{<k}\qquad \mbox{for all\,\,\,\,} -M<k\le -10\label{1.20}
\eea
 with  $U_{<k}=\sum_{k'<k}U_k$. Due to the fact that the matrices
$\tDe_k=P_k\tDe$ are antisymmetric we find the identity
$$U_k^{\,t}\cdot U_{<k}+U_{<k}^{\,t}\cdot U_k=0$$
whence,
\be{ident}
U_{<k}^{\,t}\cdot  U_{<k}-I=\sum_{k'<k}U_{k'}^{\,t}\cdot U_{k'}
\end{equation}
Using this identity we can prove inductively that
\bea
\|U_{<k}\|_{L_t^\infty L_x^\infty}&\le& 2\nn\\
\|U_{k}\|_{L_t^\infty L_x^\infty}&\lesssim& Cc_k\qquad
\mbox{for\,\,\,} k>-M \label{1.22}
\eea
as well as 
\be{1.22'}
\|U_{k}\|_{L_t^2 L_x^\infty}\lesssim C 2^{-k/2}c_k\qquad
\mbox{for\,\,\,} k>-M. 
\end{equation}
Also,
\bea
\|\pr U_{<k}\|_{L_t^\infty L_x^\infty}&\le& 2^k C^2 c_k\nn\\
\|\pr U_{k}\|_{L_t^\infty L_x^\infty}&\lesssim& 2^k C^2 c_k\label{1.23}
\eea
and 
\bea
\|\pr U_{<k}\|_{L_t^2 L_x^\infty}&\le& 2^{k/2} C^2 c_k\nn\\
\|\pr U_{k}\|_{L_t^2 L_x^\infty}&\lesssim& 2^{k/2} C^2 c_k \label{1.24}
\eea
as well as,
\bea
\|\square U_{<k}\|_{L_t^2 L_x^{n-1}}&\le& 2^{k(\frac{3k}{2}-\frac{n}{n-1})}c_k
\nn\\
\|\square U_{k}\|_{L_t^2 L_x^{n-1}}&\le& 2^{k(\frac{3k}{2}-\frac{n}{n-1})}c_k
\label{sqU} 
\eea
Indeed  the first inequality of  \eqref{1.22} holds for $k\le -M$.
Assume that it holds up to some $-M< k<-10$.  In view of part
ii) of Proposition \ref{step2} we have 
$\|\tDe_k\|_{L_t^\infty L_x^\infty}\lesssim Cc_k$.  Therefore,

$$\|U_k\|_{L_t^\infty L_x^\infty}=\|\tDe_kU_{<k}\|_{L_t^\infty L_x^\infty}\lesssim
2Cc_k
$$
which proves the second part of \eqref{1.22}. To complete the induction
for the first inequality we use the identity \eqref{ident} according to
which
$$\|U_{\le k}\|_{L_t^\infty L_x^\infty}^2\le
1+\sum_{-M<k'<k}\|U_{k'}\|_{L_t^\infty L_x^\infty}^2\le 2
$$ 
provided that $\varep$ is sufficiently small.

\qquad To prove \eqref{1.23} and \eqref{1.24}  we proceed once more by induction.
Observe  that the first estimate follows from the second.
Indeed, using \eqref{env}, 
\beaa
\|\pr U_{<k}\|_{L_t^\infty L_x^\infty}&\le &\sum_{k'<k}
\|\pr U_{k'}\|_{L_t^\infty L_x^\infty}\lesssim  C^2\sum_{k'<k}2^{k'}c_{k'}\\
&\lesssim& C^2c_k\sum_{k'<k}2^{k'}2^{\si(k-k')}=C^2c_k2^k\sum_{k''<0}
2^{k''(1-\si)}
\lesssim C^2 2^k c_k
\eeaa
as desired. 
Also,
\beaa
\|\pr U_{<k}\|_{L_t^2 L_x^\infty}&\le &\sum_{k'<k}
\|\pr U_{k'}\|_{L_t^2 L_x^\infty}\lesssim  C^2\sum_{k'<k}2^{k'/2}c_{k'}\\
&\lesssim& C^2c_k\sum_{k'<k}2^{k'/2}2^{\si(k-k')}=C^2c_k2^{k/2}\sum_{k''<0}
2^{k''(1/2-\si)}
\lesssim C^2 2^k c_k
\eeaa
since $0<\si<\frac{1}{2}$.

\qquad The second estimate in \ref{1.23} can be proved now by induction
with the help of the  definition  $U_k=\tDe_k\cdot U_{<k}$. 
The result is clearly true for $k\le -M$. We may thus assume that
the first  estimate  in \eqref{1.23} is verified 
for some $k<-10$. Using  the estimates  \eqref{1.15} and \eqref{1.16}
of part ii) of Proposition \ref{step2} we derive

\beaa
\|\pr U_k\|_{L_t^\infty L_x^\infty}&\le& \|\pr \tDe_k\|_{L_t^\infty
L_x^\infty}\|U_{<k}\|_{L_t^\infty L_x^\infty}+\| \tDe_k\|_{L_t^\infty
L_x^\infty}\|\pr U_{<k}\|_{L_t^\infty L_x^\infty}\\
&\lesssim&  2\cdot 2^kc_k +c_k\cdot  2^kC^2  c_k\lesssim 2^kc_k
\eeaa
 as desired.  For the second estimate in \eqref{1.24}
we have,
\beaa
\|\pr U_k\|_{L_t^2 L_x^\infty}&\le& \|\pr \tDe_k\|_{L_t^2
L_x^\infty}\|U_{<k}\|_{L_t^\infty L_x^\infty}+\| \tDe_k\|_{L_t^\infty
L_x^\infty}\|\pr U_{<k}\|_{L_t^2 L_x^\infty}\\
&\lesssim&  2\cdot 2^{k/2} c_k +c_k\cdot  2^{k/2} C^2  c_k\lesssim 2^kc_k
\eeaa

\qquad To prove \eqref{sqU} assume the first estimate to be true.
Then, 
$$
\square U_k=\square(\tDe_kU_{<k})=(\square\tDe_k)\cdot
U_{<k}+2\pr^\mu\tDe\cdot\pr_\mu U_{<k} 
 +\tDe_k\cdot\square U_{<k}.
$$
Hence, using the induction hypothesis and 
the estimates  we have for $U_{<k}$, $\pr U_{<k}$, $\tDe_k$ and $\pr\tDe_k$
we derive: 
\beaa
\|\square U_k\|_{L_t^2L_x^{n-1}}&\le& \|\square \tDe_k\|_{L_t^2L_x^{n-1}}
\|U_{<k}\|_{L_t^\infty L_x^\infty}\\
&+&2\|\pr\tDe_k\|_{L_t^2L_x^{n-1}}\|\pr U_{<k}\|_{L_t^\infty L_x^\infty}\\
&+&\|\tDe_k\|_{L_t^\infty L_x^\infty}\|\square U_{<k}\|_{L_t^2L_x^{n-1}}\\
&\lesssim& C 2^{k(\frac{3}{2}-\frac{n}{n-1})}c_{k}
+Cc_k2^{k(\frac 1 2-\frac{n}{n-1})}c_k\cdot 2^k Cc_k\\
 &+& C c_k\cdot C 2^{(\frac{3}{2}-\frac{n}{n-1})}c_{k}\lesssim  C
2^{k(\frac{3}{2}-\frac{n}{n-1})}c_{k}
\eeaa
Now, using \eqref{env}, 
\beaa
\|\square U_{\le k}\|_{L_t^2L_x^{n-1}}&\le&C\sum_{k'\le k}
2^{k'(\frac{3}{2}-\frac{n}{n-1})}c_{k'}\le Cc_k\sum_{k'\le k}
2^{k'(\frac{3}{2}-\frac{n}{n-1}-\si)}\\
&\lesssim& C 2^{k(\frac{3}{2}-\frac{n}{n-1})}c_k
\eeaa
as desired.

\qquad We summarize the most important properties of   $U=I+\sum_{-M<k\le-10} U_k$
in the following
\begin{proposition}
Assume that $\varep $ is sufficiently small depending on $C$
and $M$ sufficiently large depending on $T, C,\varep$. Then the matrices
$U$ verify the following properties:

\medskip
\noindent
i.)\,\,\,  Approximate  orthogonality:
\be{1.25}
\|U^{\,t}U-I\|_{L_t^\infty L_x^\infty},\,\,\,
 \|\pr (U^{\,t}U-I)\|_{L_t^\infty L_x^\infty}\lesssim C^2\varep
\end{equation}
In particular, for small $\varep$, U is invertible and we have,
\be{1.26}
\|U\|_{L_t^\infty L_x^\infty},\,\,\,
 \|U^{\,-1}\|_{L_t^\infty L_x^\infty}\lesssim 1
\end{equation}

\medskip
\noindent
ii.) \,\,\, Approximate gauge condition:
\be{1.27}
\|\pr_\mu U-\pr_\mu  \tDe\cdot U\|_{L_t^1 L_x^\infty}\lesssim C^2\varep
\end{equation}

\medskip
\noindent
iii.) \,\,\,We also have,
\be{1.28}
\|\pr U\|_{L_t^\infty L_x^\infty},\,\,\,\|\pr U\|_{L_t^1 L_x^\infty}\lesssim
C^2\varep
\end{equation}
\be{1.29}
\|\square U\|_{L_t^2 L_x^{n-1}}\lesssim C^2\varep
\end{equation}
\label{step3}
\end{proposition}
\begin{proof}:\,\,
The first part of the proposition is an easy consequence 
of the identity  \eqref{ident} as well as the estimates
\eqref{1.23}
and \eqref{1.24}. 
 To prove the  crucial second part  we write
$$\pr_\mu U-\pr_\mu\tDe\cdot U=\sum_{-M< k\le 10}
\bigg(\pr_\mu U_k-(\pr_\mu\tDe_{\le k}\cdot U_{\le k}-\pr_\mu\tDe_{< k}\cdot
U_{< k})\bigg)-\pr_\mu\tDe_{\le -M}$$
We estimate $\pr_\mu\tDe_{\le -M}$  using Cauchy-Schwartz and 
 \eqref{1.16} as follows
$$\|\pr_\mu\tDe_{\le -M}\|_{L_t^1 L_x^\infty}\le
T^{\frac{1}{2}}\|\pr_\mu\tDe_{\le -M}\|_{L_t^2 L_x^\infty}
\lesssim T^{\frac{1}{2}}\sum_{k\le -M} 2^{k/2} c_k\lesssim \varep 
T^{\frac{1}{2}} 2^{-M/2}.$$
Thus, picking $M$ sufficiently large,
$$\|\pr_\mu\tDe_{\le -M}\|_{L_t^1 L_x^\infty}\le \varep.$$

To end the proof of \eqref{1.27} it suffices  to prove that for all $-M<k\le -10$
we have
$
\|E_k\|_{L_t^1 L_x^\infty}\lesssim  C^2c_k^2
$
where
\beaa
E_k&=& \pr_\mu U_k-(\pr_\mu\tDe_{\le k}\cdot U_{\le k}-\pr_\mu\tDe_{< k}\cdot
U_{< k})\\
&=&\pr_\mu U_k-(\pr_\mu\tDe_{\le k}-\pr_\mu\tDe_{< k})\cdot U_{<k}-
\pr_\mu\tDe_{\le  k}\cdot U_k\\
&=&\pr_\mu U_k-\pr_\mu\tDe_{ k}\cdot U_{<k}-
\pr_\mu\tDe_{\le  k}\cdot U_k
\eeaa
Now using the definition $U_k=\tDe_k\cdot U_{<k}$,
\beaa
E_k&=&\pr_\mu (\tDe_k U_{<k})-(\pr_\mu\tDe_{ k}\cdot U_{<k}-
\pr_\mu\tDe_{\le  k}\cdot U_k)\\
&=&\tDe_k\cdot \pr_\mu U_{<k}+\pr_\mu\tDe_{\le  k}\cdot U_k
\eeaa
Therefore, using \eqref{1.15}, \eqref{1.16}
as well as \eqref{1.22'}, \eqref{1.24}
\beaa
\|E_k\|_{L_t^1 L_x^\infty}&\le&\|\tDe_k\|_{L_t^2 L_x^\infty}
\|\pr U_{<k}\|_{L_t^2L_x^\infty}+\|\pr \tDe_{\le k}\|_{L_t^2 L_x^\infty}
\| U_{k}\|_{L_t^2L_x^\infty}\\
&\lesssim &C^2c_k^2
\eeaa
as desired.

\qquad The   estimates  \eqref{1.28} of  part iii.)  of  the proposition
 are immediate consequences  of the estimates \eqref{1.23},
\eqref{1.24}. The inequality \eqref{1.29} can be derived immediately from
\eqref{sqU}.
\end{proof}

\qquad Following \cite{T} we  are now ready to perform the gauge transformation
\be{1.30}
\Psi=U\cdot  W
\label{gauge}
\end{equation}
W verifies the equation
\bea
\Box W&=&-2U^{-1}(\pr_\mu  U-\pr_\mu \tDe \cdot U)\pr^\mu W \label{1.31}\\
&-&2U^{-1}
\pr_\mu \tDe\cdot (\pr^\mu U) U^{-1}\Psi-U^{-1}(\square U)U^{-1}\Psi
+\mbox{error}\nn
\eea
In view of Proposition \ref{step3} we derive,

\begin{proposition} The matrix valued function $W$
verifies an equation of the form
$$\square W=\mbox{error}.$$\label{step4}
\end{proposition}
Therefore, if $\varep $ is sufficiently small, 
$$\|\Psi\|_{S_0}\lesssim \|W\|_{S_0}\le \|\Psi[0]\|_{H^{\frac{n-2}{2}}}+
CC^3\varep c_0\le  Cc_0.$$
This is precisely \eqref{1.12} which ends the proof of the Main
Proposition\footnote{ See also the more complete argument in section 7 of 
\cite{T}.}.

\begin{remark}  It is interesting to compare our  results 
 for the Hodge  system \eqref{1.2} -\eqref{1.3}  with
the system obtained by considering Lorentz gauge,  zero curvature
connections in  a  general Lie algebra, see \cite{KMa}:
 \bea
\pr_\a A_\b-\pr_\b A_\a&=&[A_\a, A_\b]\nonumber\\
\pr^\a A_\a=0.\label{zerocurv}
\eea
Such systems can be written in the form \eqref{1.2} -\eqref{1.3}
in the particular case when the structure constants $C^a_{bc}$
verify, in addition to $C^a_{bc}=-C^a_{cb}$, the relations
$C^a_{bc}=-C^b_{ac}$. In this case $\Ga^a_{bc}=C^a_{bc}$. This corresponds
to the case of a Lie group with a bi-invariant  Riemannian metric such as 
$S^3$. The 
system  \eqref{zerocurv} is interesting however  in its own right,
for general Lie algebras. The results of this paper can be extended
to the case of classical  Lie algebras such as $o(n), su(n)$
 and probably more generally to Lie algebras of
 compact Lie groups. The compactness
seems in this case to be essential\footnote{In this case the  transformation
\eqref{gauge} should be replaced by the partial gauge transformation 
$A\longrightarrow UAU^{-1}$ with $U$ an element of the group. Compactness
of the group is neede to control  
 the sup-norm of $U$.}, by contrast to the case
of wave  maps where the compactness of the target manifold is not important.
\end{remark}

\end{document}